\pdfoutput=1
\documentclass[11pt]{amsart}
\usepackage{amsmath,amssymb,amsthm}
\usepackage[margin=1in]{geometry}
\usepackage{hyperref}
\usepackage{booktabs}
\usepackage{placeins}

\newtheorem{theorem}{Theorem}
\newtheorem{lemma}[theorem]{Lemma}
\newtheorem{proposition}[theorem]{Proposition}
\newtheorem{corollary}[theorem]{Corollary}
\theoremstyle{definition}

\theoremstyle{remark}
\newtheorem{remark}[theorem]{Remark}

\newcommand{\fl}[1]{\lfloor #1\rfloor}
\newcommand{\fr}[1]{\{#1\}}
\newcommand{\one}{\mathbf 1}
\hypersetup{pdftitle={Discrepancy of the words floor((2n+k)r)-floor((n+k)r)-floor(nr) and twenty-two conjectures of Kimberling}, pdfauthor={Alex Ashburn}}

\title[Kimberling's $2n-a(n)$ conjectures]{Discrepancy of the words $\fl{(2n+k)r}-\fl{(n+k)r}-\fl{nr}$\\ and twenty-two conjectures of Kimberling}

\author{Alex Ashburn}
\address{Independent Researcher}
\email{aashburn1998@gmail.com}

\subjclass[2020]{11K38, 11J70, 11B85, 68R15}
\keywords{irrational rotation, local discrepancy, one-sided boundedness, continued fraction, Ostrowski numeration, automatic sequence, Walnut, OEIS}

\begin{document}

\begin{abstract}
For an irrational $r>0$ and an integer $k\ge1$, the word
$s_k(n)=\fl{(2n+k)r}-\fl{(n+k)r}-\fl{nr}$ takes only the values $0$ and $1$, each with density $1/2$.
In 2019 Kimberling added to the OEIS eleven such words, for $r\in\{\varphi,\sqrt2,\sqrt3,e\}$ and $k\le4$, together with the $22$ sequences $a(n)$ listing the positions of their $0$s and $1$s.
For $21$ of these he conjectured that $2n-a(n)$ is unbounded below and above, and for the remaining one that $a(n)-2n$ is unbounded above.
We settle all $22$ conjectures: $19$ are true and $3$ are false.
The false ones all come from $r=\sqrt2$ with $k\in\{2,4\}$: the partial sums of $2s_k(n)-1$ never exceed $0$ (for $k=2$) or $1$ (for $k=4$).
Hence $a(n)\ge2n$ for A327207, $2n-a(n)\ge1$ for A327206, and $2n-a(n)\ge0$ for A327223.
For even $k$ the qualitative statements follow, through a short reduction, from a theorem of Ying and Zheng on the one-sided boundedness of the local discrepancy $D_N(\alpha,1/2)$, and also from earlier formulas of Ro\c{c}adas and Schoi\ss{}engeier, which Amoroso and Omarjee applied to $r=e$; for $r\in\{\varphi,\sqrt2,\sqrt3\}$ they also follow from results of Boshernitzan and Ralston and of Dolgopyat and Sarig.
The cases with odd $k$, where that reduction leads to a sum of two local discrepancies with irrational cut points, are new, and so are the exact bounds.
For quadratic $r$ we prove the results with Ostrowski-automatic sequences, verified in Walnut, and an exact block-sum analysis that also gives the growth rates of the extremes; for even $k$ these rates also follow from known results on $D_N(\alpha,1/2)$.
For $r=e$ we prove an exact formula for the discrepancy at $N=tq_i$, where $q_i$ are the convergent denominators of $e$, and obtain closed forms along $N=tq_{3j-2}$.
\end{abstract}

\maketitle

\section{Introduction}

Fix an irrational number $r>0$ and an integer $k\ge1$, and put
\begin{equation}\label{eq:s}
  s_k(n)\;=\;\fl{(2n+k)r}-\fl{(n+k)r}-\fl{nr}\qquad(n\ge0).
\end{equation}
Lemma~\ref{lem:rot} below shows that $s_k(n)\in\{0,1\}$ and that both letters have density $1/2$.
In August and September 2019 Kimberling entered into the OEIS~\cite{oeis} the eleven words listed in Table~\ref{tab:words}, and for each of them the sequences of positions $n\ge1$ of its $0$s and of its $1$s; Jianing Song corrected the position sequences in September 2019.
Each of these $22$ position sequences $a(n)$ satisfies $a(n)/n\to2$, and the entry for each one states:
\begin{quote}
\emph{Although $a(n)/n\to2$, the sequence $2n-a(n)$ appears to be unbounded below and above.}
\end{quote}
The one exception is A327224, whose comment says only that $a(n)-2n$ appears to be unbounded above.
As of September 2026 these conjectures are still marked as open in the OEIS.

\begin{table}[ht]
\centering
\begin{tabular}{llcll}
\toprule
word & $r$ & $k$ & positions of $0$ & positions of $1$\\
\midrule
A327174 & $\varphi=(1+\sqrt5)/2$ & 1 & A327175 & A327176\\
A327177 & $\sqrt2$ & 1 & A327178 & A327179\\
A327180 & $\sqrt3$ & 1 & A327181 & A327182\\
A327183 & $e$ & 1 & A327184 & A327185\\
A327202 & $\varphi$ & 2 & A327203 & A327204\\
A327205 & $\sqrt2$ & 2 & A327206 & A327207\\
A327208 & $\sqrt3$ & 2 & A327209 & A327210\\
A327211 & $e$ & 2 & A327212 & A327213\\
A327216 & $\varphi$ & 3 & A327217 & A327218\\
A327219 & $\varphi$ & 4 & A327220 & A327221\\
A327222 & $\sqrt2$ & 4 & A327223 & A327224\\
\bottomrule
\end{tabular}
\caption{Kimberling's words $s_k$ (offset $0$) and their position sequences (offset $1$).}
\label{tab:words}
\end{table}

Let
\[
  D_k(N)\;=\;\sum_{n=1}^N\bigl(2s_k(n)-1\bigr)\qquad(N\ge0)
\]
be the excess of $1$s over $0$s among $s_k(1),\dots,s_k(N)$.
Lemma~\ref{lem:equiv} shows that each of the $21$ two-sided conjectures is equivalent to the statement that $D_k$ is unbounded above and unbounded below.
The conjecture for A327224 is equivalent to $D_k$ being unbounded below.

\begin{theorem}\label{thm:main}
\begin{enumerate}
\item For $(r,k)\in\{(\varphi,1),(\varphi,2),(\varphi,3),(\varphi,4),(\sqrt2,1),(\sqrt3,1),(\sqrt3,2),(e,1),(e,2)\}$ the sequence $D_k(N)$ is unbounded above and unbounded below.
\item For $r=\sqrt2$ and $k\in\{2,4\}$ the sequence $D_k(N)$ is unbounded below, and
\[
  \sup_{N\ge0}D_2(N)=0,\qquad \sup_{N\ge0}D_4(N)=1 .
\]
\end{enumerate}
Consequently, $19$ of the $22$ conjectures are true.
The conjectures for A327206, A327207 and A327223 are false: for all $n\ge1$,
\[
  2n-a(n)\ge1\ \text{(A327206)},\qquad 2n-a(n)\le0\ \text{(A327207)},\qquad 2n-a(n)\ge0\ \text{(A327223)},
\]
and each bound is attained, for instance at $n=1$, $n=1$ and $n=2$ respectively.
For A327224 the conjecture is true, and moreover $a(n)-2n\ge-1$ for all $n$.
\end{theorem}

\subsection*{Relation to earlier work}
By Lemma~\ref{lem:rot}, $D_k(N)$ is twice the discrepancy of the rotation $n\mapsto\fr{nr}$ with respect to a union of two intervals of total length $1/2$.
For a single interval $[0,c)$ the local discrepancy is
\[
  D_N(\alpha,c)\;=\;\#\{1\le n\le N:\fr{n\alpha}<c\}-cN .
\]
By the classical results of Hecke~\cite{hecke} and Ostrowski~\cite{ostrowski}, $D_N(\alpha,c)$ is bounded when $c\in\mathbb Z\alpha+\mathbb Z$, and Kesten~\cite{kesten} proved the converse.
One-sided boundedness is subtler.
It was studied by Sós~\cite{sos}, by Dupain~\cite{dupain}, and by Dupain and Sós~\cite{dupain-sos}.
For rational $c$, Ying and Zheng~\cite{yz} gave a complete criterion in terms of the convergents of $\alpha$ (Theorem~\ref{thm:yz} below).
The symmetric case $c=1/2$ has been studied in more detail.
Boshernitzan and Ralston~\cite{br} proved that $D_N(\alpha,1/2)\ge0$ for all $N$ if and only if every partial quotient $a_j$ of $\alpha$ with odd $j$ is even.
Ro\c{c}adas and Schoi\ss{}engeier~\cite{rs} gave, for every $c$ and up to a bounded error, the maximum and the minimum of $D_N(\alpha,c)$ over $1\le N<q_{m+1}$; for $c=1/2$ these are $\tfrac14$ times the sum of the partial quotients $a_{j+1}$, $j\le m$, taken over the even $j$ with $q_j$ odd for the maximum and over the odd $j$ with $q_j$ odd for the minimum (with a minus sign).
Amoroso and Omarjee~\cite{ao} applied these formulas to $\alpha=e$ and $\alpha=e/2$; in particular, their proof shows that $D_N(e,1/2)$ is unbounded above and below, with extremes of order $(\log N/\log\log N)^2$.
O'Bryant, Reznick and Serbinowska~\cite{ors} studied $S_N(\beta)=\sum_{n=1}^N(-1)^{\fl{n\beta}}$, which equals $2D_N(\beta/2,1/2)$.
For $\beta=[b_0;b_1,\dots]$ with $b_0$ odd and all other $b_i$ even, they gave an exact formula for $S_N(\beta)$ in terms of an Ostrowski-type expansion of $N$ and the least $N$ with $S_N(\beta)=j$ for every $j$~\cite[Theorem~4.4 and Corollary~4.5]{ors}; this covers $2\varphi=[3;\overline4]$ and $2\sqrt3=[3;\overline{2,6}]$.
Bruin and Fokkink~\cite{bf} answered two questions of Kimberling on the OEIS entry A120243, which lists the $n$ with $\fr{n\sqrt2}<1/2$ (its complement A120749 lists those with $\fr{n\sqrt2}>1/2$); the first of them follows from the criterion of~\cite{br}.
They also showed~\cite[Theorem~5]{bf} that the records of $S_N(2\sqrt2)=2D_N(\sqrt2,1/2)$ occur exactly at the half-Pell numbers $1,6,35,204,1189,\dots$.
Ralston~\cite{ralston} studied, for fixed $\alpha$, the set of starting points $x$ for which the analogous sums along the orbit of $x$ are never negative, and Dolgopyat and Sarig~\cite{ds} recently determined, for $\alpha$ of bounded type, when the values of $D_N(\alpha,1/2)$ are equidistributed in $\tfrac12\mathbb Z$.

For even $k$ we show in Proposition~\ref{prop:even} that $D_k(N)=-2D_N(r,1/2)+O(1)$.
Combined with the Ying--Zheng criterion, this already decides all conjectures with even $k$, except for the exact values of the suprema in part~(2).
We record this in Corollary~\ref{cor:yz} and credit these qualitative statements to~\cite{yz}; for $c=1/2$, which is all we need, they can also be derived from the earlier formulas of~\cite{rs}, and for $r=e$ from their application in~\cite{ao}.
For $r\in\{\varphi,\sqrt2,\sqrt3\}$ they also follow from~\cite{br} and~\cite{ds} (Remark~\ref{rem:br-ds}).
For odd $k$ the same reduction produces the irrational cut points $c_1=(1-\fr{kr})/2$ and $c_1+\tfrac12$, which are congruent to $\fr{-r/2}$ and $\fr{1/2-r/2}$ modulo $\mathbb Zr+\mathbb Z$.
They enter~\eqref{eq:decomp} with the same sign, so $D_k$ is, up to a bounded error, minus twice the sum of two local discrepancies, not the discrepancy of the translate $[c_1,c_1+\tfrac12)$ of $[0,1/2)$.
Equivalently, the identity $\one[x<c_1]+\one[x<c_1+\tfrac12]=\one[x<\tfrac12]+\one[\fr{2x}<2c_1]$ for $x\in[0,1)$ gives
\[
  D_k(N)=-2D_N(r,\tfrac12)-2D_N(2r,1-\fr{kr})+2D_N(r,1-\fr{kr}),
\]
where the last term is bounded by Lemma~\ref{lem:arc}.
For odd $k$ the point $1-\fr{kr}$ is not in $2r\mathbb Z+\mathbb Z$, so by Kesten's theorem the middle term is unbounded.
Theorem~\ref{thm:yz} does not cover this case.
Neither do the results of~\cite{br,ds,ors,bf}, which concern intervals $[0,c)$ with rational $c$, or those of~\cite{ralston}, which concern translates of $[0,1/2)$; the formulas of~\cite{rs} hold for every $c$, but they describe the extremes of a single local discrepancy, not of a sum of two of them as above.
So we use other methods.
\begin{itemize}
\item For $r\in\{\varphi,\sqrt2,\sqrt3\}$ (Section~\ref{sec:quad}) the words $s_k$ are Ostrowski-automatic.
We construct deterministic finite automata with output (DFAOs) and prove them correct with the Walnut prover~\cite{walnut,shallit-book}, using the method of Schaeffer, Shallit and Zorcic~\cite{ssz} for $\fl{n\gamma}$.
We then compute the partial sums $D_k$ exactly through a finite linear recursion for block sums.
This gives explicit families of $N$ along which $D_k$ grows linearly in $\log N$, the exact suprema in part~(2), and the exact growth rates of $\max_{M\le X}D_k(M)$ and $\min_{M\le X}D_k(M)$ (Proposition~\ref{prop:growth}); for even $k$ these rates also follow from~\cite{rs,ors,bf} (Remark~\ref{rem:rates}).
\item For $r=e$ (Section~\ref{sec:e}) we prove an exact formula for $D_N(\alpha,c)$ when $N=tq_i$ with $t\le a_{i+1}$ (Lemma~\ref{lem:block}).
From it we derive closed forms for $D_1$ and $D_2$ along $N=tq_{3j-2}$ (Theorem~\ref{thm:e}).
\end{itemize}
Everything was checked independently by computer; see Section~\ref{sec:verification}.

\section{Preliminaries}\label{sec:prelim}

\subsection{Positions and partial sums}

\begin{lemma}\label{lem:equiv}
Let $s:\mathbb Z_{\ge1}\to\{0,1\}$ take both values infinitely often, let $D(N)=\sum_{n\le N}(2s(n)-1)$, and let $u(m)$ and $v(m)$ be the positions of the $m$-th $1$ and the $m$-th $0$.
Then
\[
  2m-u(m)=D(u(m)),\qquad 2m-v(m)=-D(v(m)).
\]
Moreover, $2m-u(m)$ is unbounded above (below) if and only if $D$ is unbounded above (below), and $2m-v(m)$ is unbounded above (below) if and only if $D$ is unbounded below (above).
\end{lemma}

\begin{proof}
Among $s(1),\dots,s(u(m))$ there are $m$ ones and $u(m)-m$ zeros, so $D(u(m))=2m-u(m)$.
The identity for $v$ is similar.
If $D(N)>\max_{n<N}D(n)$ then $s(N)=1$, so every new maximum of $D$ is attained at a position of a $1$; and $D(u(m))\le\sup D$.
If $D(N)=-M$, let $u$ be the first position of a $1$ after $N$; then $D(u)=D(u-1)+1\le D(N)+1$.
Hence $\inf_mD(u(m))\le\inf D+1$.
The statements for $v$ follow by exchanging the roles of $0$ and $1$.
\end{proof}

\subsection{Rotation form}

Write $\fr{x}=x-\fl x$.
Since $\fl{mr}=m\fl r+\fl{m\fr r}$, the word $s_k$ depends only on $\gamma=\fr r$, and $\fr{nr}=\fr{n\gamma}$.

\begin{lemma}\label{lem:rot}
Let $\rho=\fr{kr}$ and
\[
  c_1=\frac{1-\rho}{2},\qquad c_2=1-\rho,\qquad c_3=1-\frac{\rho}{2}.
\]
For $n\ge1$ we have $s_k(n)=\one_{[c_1,c_2)\cup[c_3,1)}(\fr{nr})$, and for $N\ge0$
\begin{equation}\label{eq:decomp}
  D_k(N)=2\bigl(D_N(r,c_2)-D_N(r,c_1)-D_N(r,c_3)\bigr).
\end{equation}
\end{lemma}

\begin{proof}
Put $x=\fr{nr}$.
Then $\fl{(2n+k)r}=2\fl{nr}+\fl{kr}+\fl{2x+\rho}$ and $\fl{(n+k)r}=\fl{nr}+\fl{kr}+\fl{x+\rho}$, so $s_k(n)=\fl{2x+\rho}-\fl{x+\rho}$.
If $x+\rho<1$, this equals $\one[2x+\rho\ge1]=\one[x\ge c_1]$, and $x<c_2$.
If $x+\rho\ge1$, that is $x\ge c_2$, it equals $\one[2x+\rho\ge2]=\one[x\ge c_3]$.
Since $c_1<c_2<c_3$, this proves the first claim.
(For $n\ge1$ the point $\fr{nr}$ never equals a $c_i$: for example, $x=c_1$ would give $(2n+k)r\in\mathbb Z$.)
For the second claim, write $C(c)=\#\{n\le N:\fr{nr}<c\}=D_N(r,c)+cN$.
Then $\sum_{n\le N}s_k(n)=C(c_2)-C(c_1)+N-C(c_3)$.
The coefficient of $N$ in $2\sum s_k(n)-N$ is $2(c_2-c_1+1-c_3)-1=0$.
\end{proof}

In particular $s_k(n)\in\{0,1\}$, and the set $[c_1,c_2)\cup[c_3,1)$ has measure $1/2$.

\subsection{Continued fractions}\label{sec:cf}

Let $\alpha=[a_0;a_1,a_2,\dots]$ be irrational, with convergents $p_i/q_i$ ($p_{-1}=1$, $q_{-1}=0$, $p_0=a_0$, $q_0=1$), complete quotients $\alpha_i=[a_i;a_{i+1},\dots]$, and $\theta_i=q_i\alpha-p_i$.
We use the standard facts~\cite{rockett-szusz}:
\begin{align}
  \theta_i&=\frac{(-1)^i}{q_i\alpha_{i+1}+q_{i-1}},\quad\text{so}\quad \frac{1}{q_i|\theta_i|}=\alpha_{i+1}+\frac{q_{i-1}}{q_i};\label{eq:theta}\\
  p_iq_{i-1}-p_{i-1}q_i&=(-1)^{i-1};\label{eq:det}\\
  |\theta_{i-1}|&=a_{i+1}|\theta_i|+|\theta_{i+1}|,\qquad |\theta_{i-1}|<1/q_i\quad(i\ge1).\label{eq:rec}
\end{align}

\begin{lemma}[Hecke, Ostrowski]\label{lem:arc}
Let $m\ge1$ be an integer and $\beta$ real.
Then for every $N\ge1$,
\[
  \Bigl|\#\{1\le n\le N:\fr{n\alpha+\beta}\ge1-\fr{m\alpha}\}-N\fr{m\alpha}\Bigr|<m .
\]
\end{lemma}

\begin{proof}
For real $x,y$ we have $\one[\fr x+\fr y\ge1]=\fr x+\fr y-\fr{x+y}$.
Taking $x=n\alpha+\beta$, $y=m\alpha$ and summing over $n\le N$, the count equals
$N\fr{m\alpha}+\sum_{n=1}^{m}\fr{n\alpha+\beta}-\sum_{n=N+1}^{N+m}\fr{n\alpha+\beta}$.
Each of the two sums lies in $[0,m)$.
\end{proof}

\section{Even \texorpdfstring{$k$}{k} and the theorem of Ying and Zheng}\label{sec:even}

\begin{proposition}\label{prop:even}
Let $k=2l$ be even.
Then $D_k(N)=-2D_N(r,1/2)+R_k(N)$ with $|R_k(N)|<4k$ for all $N\ge1$.
\end{proposition}

\begin{proof}
Let $\sigma=\fr{lr}$ and $\varepsilon=\fl{2\sigma}\in\{0,1\}$, so that $\rho=\fr{kr}=2\sigma-\varepsilon$.
First, $c_2=1-\rho$, and $D_N(r,c_2)=-(A-N\rho)$ with $A=\#\{n\le N:\fr{nr}\ge1-\rho\}$.
Lemma~\ref{lem:arc} with $m=k$, $\beta=0$ gives $|D_N(r,c_2)|<k$.
Next, $c_1=\tfrac12-\sigma+\tfrac\varepsilon2$ and $c_3=1-\sigma+\tfrac\varepsilon2$.
So $\{c_1,c_3\}=\{1-\sigma,\,t\}$, where $t=\tfrac12-\sigma$ if $\varepsilon=0$ and $t=\tfrac32-\sigma$ if $\varepsilon=1$.
As before, $|D_N(r,1-\sigma)|<l$.
If $\varepsilon=0$, then $t\le\fr{nr}<\tfrac12$ is equivalent to $\fr{nr+\tfrac12}\ge1-\sigma$.
Hence $D_N(r,\tfrac12)-D_N(r,t)=\#\{n:\fr{nr+\tfrac12}\ge1-\sigma\}-N\sigma$, which has absolute value less than $l$ by Lemma~\ref{lem:arc} (with $m=l$, $\beta=\tfrac12$).
If $\varepsilon=1$, then $\tfrac12\le\fr{nr}<t$ is equivalent to $\fr{nr+\tfrac12}<1-\sigma$.
Hence $D_N(r,t)-D_N(r,\tfrac12)=-(\#\{n:\fr{nr+\tfrac12}\ge1-\sigma\}-N\sigma)$, which is again smaller than $l$ in absolute value.
By \eqref{eq:decomp}, $D_k(N)=-2D_N(r,\tfrac12)+R_k(N)$ with
$R_k=2D_N(r,c_2)-2D_N(r,1-\sigma)-2\bigl(D_N(r,t)-D_N(r,\tfrac12)\bigr)$.
Therefore $|R_k|<2k+2l+2l=4k$.
\end{proof}

Numerically $|R_k(N)|\le2$ for all four $r$, $k\in\{2,4\}$ and $N\le10^6$.

\begin{theorem}[Ying and Zheng~{\cite[Theorem~1.2]{yz}}]\label{thm:yz}
Let $\alpha>0$ be irrational with convergent denominators $q_n$ (with $q_{-1}=0$), and let $c=h/k'\in(0,1)$ be rational in lowest terms.
Then $D_N(\alpha,c)$ is bounded above (respectively below) if and only if there is an even (respectively odd) $m\ge-1$ such that $k'\mid q_{m+2n}$ for all $n\ge0$.
\end{theorem}

\begin{corollary}\label{cor:yz}
Let $k$ be even.
For $r\in\{\varphi,\sqrt3,e\}$, $D_k$ is unbounded above and below.
For $r=\sqrt2$, $D_k$ is bounded above and unbounded below.
\end{corollary}

\begin{proof}
Apply Theorem~\ref{thm:yz} with $c=1/2$ and Proposition~\ref{prop:even}: $D_k$ is bounded above if and only if $D_N(r,\tfrac12)$ is bounded below, which happens if and only if all $q_{m+2n}$ ($n\ge0$) are even for some odd $m\ge-1$.
The same holds with ``below'', ``above'' and ``even $m$''.
For $\varphi$, $q_n=F_{n+1}$ is even if and only if $3\mid n+1$, so no progression with difference $2$ consists of even terms.
For $\sqrt3=[1;\overline{1,2}]$, $q_n$ is even exactly when $n\equiv3\pmod4$ (including $q_{-1}=0$).
For $e$, $q_n$ is even exactly when $n\equiv3,5\pmod 6$ (Section~\ref{sec:e}).
For $\sqrt2=[1;\overline{2}]$, every $q_{2n-1}$ is even and every $q_{2n}$ is odd; so $m=-1$ works for ``bounded below'' and no even $m$ works.
\end{proof}

\begin{remark}\label{rem:br-ds}
For $r\in\{\varphi,\sqrt2,\sqrt3,e\}$, Corollary~\ref{cor:yz} also follows from earlier and independent work.
Since $\sqrt2=[1;\overline2]$, the criterion of Boshernitzan and Ralston~\cite[Theorem~1]{br} gives $D_N(\sqrt2,1/2)\ge0$ for all $N$, so $D_k$ is bounded above for even $k$ by Proposition~\ref{prop:even}.
In the same way, Bruin and Fokkink~\cite[Corollary~2]{bf} deduced from~\cite{br} that $a(n)<b(n)$ for all $n$, where $a$ and $b$ are A120243 and A120749.
For even $k$, $D_k$ is unbounded below, because $D_N(\sqrt2,1/2)$ is unbounded by Kesten's theorem~\cite{kesten}.
For $\alpha$ of bounded type, Dolgopyat and Sarig~\cite[Theorem~1.2]{ds} show that the values of $\#\{0\le n<N:\fr{n\alpha}<1/2\}-N/2$ are equidistributed in $\frac12\mathbb Z$ if and only if a condition on the signs of the partial quotients of the even continued fraction expansion of $2\alpha-1$ holds.
The condition holds for $\alpha=\sqrt3$, which is an example in~\cite{ds}, and for $\alpha=\varphi$, because $2\varphi-1=\sqrt5=[2;\overline4]$ has only positive even partial quotients.
Equidistribution implies that every value in $\frac12\mathbb Z$ is taken, so for these two numbers $D_N(r,1/2)$, and hence $D_k$ for even $k$, is unbounded in both directions.
For $r=e$, the computation of Amoroso and Omarjee~\cite[proof of Theorem~1.3]{ao}, based on~\cite{rs}, shows that $D_N(e,1/2)$ is unbounded in both directions, so $D_k$ is too for even $k$.
\end{remark}

Corollary~\ref{cor:yz} proves the conjectures for A327203, A327204, A327209, A327210, A327212, A327213, A327220, A327221 and A327224.
It disproves those for A327206, A327207 and A327223.
Our own proofs in Sections~\ref{sec:e} and~\ref{sec:quad} do not use~\cite{yz}.
They also give the exact suprema in Theorem~\ref{thm:main}(2).

\section{A block formula, and the case \texorpdfstring{$r=e$}{r = e}}\label{sec:e}

\begin{lemma}\label{lem:block}
Let $\alpha$ be irrational, $i\ge1$, and write $q=q_i$, $p=p_i$, $\theta=\theta_i$.
Let $1\le t\le a_{i+1}$, and let $c\in(0,1)$ with $qc\notin\mathbb Z$.
For $G\in\mathbb Z$ let $m(G)\in\{1,\dots,q\}$ be the solution of $m\,p\equiv G\pmod q$.
\begin{enumerate}
\item If $\theta>0$, put $G=\fl{qc}$ and $\lambda=\dfrac{\fr{qc}}{q\theta}-\dfrac{m(G)}{q}$. Then
\[
  D_{tq}(\alpha,c)=\min\bigl(t,\max(0,\lceil\lambda\rceil)\bigr)-t\fr{qc}.
\]
\item If $\theta<0$, put $G=\lceil qc\rceil$ and $\mu=\dfrac{\fr{-qc}}{q|\theta|}-\dfrac{m(G)}{q}$. Then
\[
  D_{tq}(\alpha,c)=t\fr{-qc}-\min\bigl(t,\max(0,\fl\mu+1)\bigr).
\]
\end{enumerate}
\end{lemma}

\begin{proof}
Write $n=m+sq$ with $1\le m\le q$ and $0\le s\le t-1$; this lists $n=1,\dots,tq$ once each.
Since $q\alpha=p+\theta$ and $m\alpha=mp/q+m\theta/q$,
\[
  \fr{n\alpha}=\Bigl\{\frac{g_m}{q}+\Bigl(\frac mq+s\Bigr)\theta\Bigr\},\qquad g_m\equiv mp\pmod q,\ 0\le g_m<q .
\]
The map $m\mapsto g_m$ is a bijection onto $\{0,\dots,q-1\}$ with $g_q=0$, because $\gcd(p,q)=1$.
The perturbation satisfies $0<(m/q+s)|\theta|\le t|\theta|\le a_{i+1}|\theta_i|<|\theta_{i-1}|<1/q$ by~\eqref{eq:rec}.

(1) If $\theta>0$, then $\fr{n\alpha}\in(g_m/q,(g_m+1)/q)$.
Such a point lies below $c\in[G/q,(G+1)/q)$ if $g_m<G$, and does not if $g_m>G$.
For $g_m=G$ (that is, $m=m(G)$) it lies below $c$ if and only if $(m/q+s)\theta<\fr{qc}/q$, that is, $s<\lambda$.
So the count is $tG+\min(t,\max(0,\lceil\lambda\rceil))$, and subtracting $tqc=tG+t\fr{qc}$ gives the claim.

(2) If $\theta<0$, then $\fr{n\alpha}\in((g_m-1)/q,g_m/q)$ for $g_m\ge1$, and $\fr{n\alpha}\in(1-1/q,1)$ for $g_m=0$.
Now $c\in((G-1)/q,G/q)$ with $1\le G\le q$.
The points with $1\le g_m\le G-1$ lie below $c$.
The points with $G+1\le g_m\le q-1$ do not, and neither does $g_m=0$ when $G<q$.
For $g_m\equiv G\pmod q$ the point equals $G/q-(m/q+s)|\theta|$, which lies below $c$ if and only if $(m/q+s)|\theta|>\fr{-qc}/q$, that is, $s>\mu$.
So the count is $t(G-1)+t-\min(t,\max(0,\fl\mu+1))$, and subtracting $tqc=tG-t\fr{-qc}$ gives the claim.
\end{proof}

We now take $\alpha=e$.
By Euler's expansion (see~\cite{cohn} for a short proof),
\[
  e=[2;1,2,1,1,4,1,1,6,\dots],\qquad a_{3l+1}=a_{3l+3}=1,\quad a_{3l+2}=2l+2\quad(l\ge0).
\]

\begin{theorem}\label{thm:e}
Let $r=e$, $j\ge2$, $i=3j-2$ and $1\le t\le 2j$ $(=a_{i+1})$.
Then
\[
  D_1(tq_i)=\begin{cases} t-2\min(t,j), & j\text{ even},\\ 2\min(t,j+1)-t, & j\text{ odd},\end{cases}
  \qquad
  D_2(tq_i)=\begin{cases} t-2\min(t,j+1), & j\text{ even},\\ 2\min(t,j)-t, & j\text{ odd}.\end{cases}
\]
In particular $D_1(jq_{3j-2})=-j$ and $D_2((j+1)q_{3j-2})=-(j+1)$ for even $j$, while $D_1((j+1)q_{3j-2})=j+1$ and $D_2(jq_{3j-2})=j$ for odd $j$.
Hence $D_1$ and $D_2$ are unbounded above and below.
\end{theorem}

For example, $q_4=7$, $q_7=71$ and $q_{10}=1001$ give $D_1(14)=-2$, $D_1(284)=4$, $D_1(4004)=-4$, $D_2(21)=-3$, $D_2(213)=3$ and $D_2(5005)=-5$.

\begin{proof}
Write $q=q_i$, $p=p_i$, $\theta=\theta_i$, $q'=q_{i-1}$ and $\alpha'=\alpha_{i+1}$.
We need five facts.

\emph{(E1) $2j+\tfrac12<\alpha'<2j+\tfrac23$.}
We have $\alpha'=[2j;1,1,\alpha_{i+4}]$ with $\alpha_{i+4}>1$.
So $\alpha_{i+3}=1+1/\alpha_{i+4}\in(1,2)$, then $\alpha_{i+2}=1+1/\alpha_{i+3}\in(\tfrac32,2)$, and finally $\alpha'=2j+1/\alpha_{i+2}$.

\emph{(E2) Parities.}
Since $a_n\bmod2$ is periodic with period $(1,0,1)$ for $n\ge1$, induction gives that $q_n$ is even exactly for $n\equiv3,5\pmod6$ and $p_n$ is even exactly for $n\equiv0,2\pmod6$ ($n\ge0$).
Hence $q$ and $p$ are odd, and $q'=q_{3j-3}$ is even if and only if $j$ is even.

\emph{(E3)} $\theta$ has the sign $(-1)^i=(-1)^j$, and $1/(q|\theta|)=\alpha'+q'/q$ by~\eqref{eq:theta}.
Also $|\theta|<1/q_{i+1}\le\tfrac12$.

\emph{(E4)} By~\eqref{eq:det}, $pq'\equiv(-1)^{i-1}=(-1)^{j-1}\pmod q$, so $p^{-1}\equiv(-1)^{j-1}q'\pmod q$.

\emph{(E5)} $q'\ge4$, and $0\le q'-\tfrac{q+1}{2}<q$.
Indeed $a_i=a_{i-1}=1$ give $q=q'+q_{i-2}$ and $q'=q_{i-2}+q_{i-3}$, so $2q'-q-1=q_{i-3}-1\ge0$.

Put $x_1=3-e$, $x_2=\tfrac{3-e}{2}$, $x_3=2-\tfrac e2$, $x_4=6-2e$ and $x_5=\tfrac72-e$.
By Lemma~\ref{lem:rot} (with $\rho=e-2$ for $k=1$ and $\rho=2e-5$ for $k=2$),
\begin{equation}\label{eq:e-decomp}
  D_1=2\bigl(D(x_1)-D(x_2)-D(x_3)\bigr),\qquad D_2=2\bigl(D(x_4)-D(x_1)-D(x_5)\bigr),
\end{equation}
where $D(x)=D_{tq}(e,x)$.
We evaluate each $D(x_\nu)$ with Lemma~\ref{lem:block}.
Throughout, $p^{-1}$ is taken modulo $q$ and we use (E3) to simplify $\lambda$ and $\mu$.

\smallskip\noindent\emph{Case $j$ even ($\theta>0$, $p^{-1}\equiv-q'$, $q'$ even).}
\begin{itemize}
\item $x_1$: $qx_1=(3q-p)-\theta$, so $\fr{qx_1}=1-\theta$ and $G\equiv-p-1$.
Then $m(G)\equiv-1-p^{-1}\equiv q'-1$, so $m(G)=q'-1$.
This gives $\lambda=\alpha'$, $\lceil\lambda\rceil=2j+1>t$, and $D(x_1)=t\theta$.
\item $x_4$: $qx_4=(6q-2p)-2\theta$, so $\fr{qx_4}=1-2\theta$ and $m(G)\equiv-2-p^{-1}\equiv q'-2$.
This gives $\lambda=\alpha'$ and $D(x_4)=2t\theta$.
\item $x_2$: $qx_2=\tfrac{3q-p}{2}-\tfrac\theta2$, so $\fr{qx_2}=1-\tfrac\theta2$ and $2G\equiv-p-2$.
Then $2m(G)\equiv-1-2p^{-1}\equiv2q'-1$, so $m(G)\equiv q'-\tfrac{q+1}2$.
By (E5), $m(G)$ equals $q'-\tfrac{q+1}2$ or $q'-\tfrac{q+1}2+q$.
Correspondingly $\lambda=\alpha'+\tfrac12$ or $\alpha'-\tfrac12$; in both cases $\lambda>2j\ge t$, so $D(x_2)=t-t(1-\tfrac\theta2)=\tfrac{t\theta}2$.
\item $x_3$: $qx_3=\bigl(2q-\tfrac{p+1}2\bigr)+\tfrac{1-\theta}2$, so $\fr{qx_3}=\tfrac{1-\theta}2$ and $2G\equiv-p-1$.
Then $2m(G)\equiv-1-p^{-1}\equiv q'-1$.
Since $q'-1$ is odd, $m(G)=\tfrac{q+q'-1}2$.
Hence $\lambda=\tfrac{\alpha'-1}2\in(j-\tfrac14,j-\tfrac16)$, $\lceil\lambda\rceil=j$, and $D(x_3)=\min(t,j)-t\tfrac{1-\theta}2$.
\item $x_5$: $qx_5=\bigl(\tfrac{7q-1}2-p\bigr)+\tfrac12-\theta$, so $\fr{qx_5}=\tfrac12-\theta$ and $2G\equiv-1-2p$.
Then $2m(G)\equiv-p^{-1}-2\equiv q'-2$, so $m(G)=\tfrac{q'}2-1\ge1$.
Hence $\lambda=\tfrac{\alpha'}2\in(j+\tfrac14,j+\tfrac13)$, $\lceil\lambda\rceil=j+1$, and $D(x_5)=\min(t,j+1)-t(\tfrac12-\theta)$.
\end{itemize}

\smallskip\noindent\emph{Case $j$ odd ($\theta<0$, $p^{-1}\equiv q'$, $q'$ odd).}
\begin{itemize}
\item $x_1$: $qx_1=(3q-p)+|\theta|$, so $G=3q-p+1$ and $\fr{-qx_1}=1-|\theta|$.
Then $m(G)\equiv p^{-1}-1\equiv q'-1$, so $\mu=\alpha'$, $\fl\mu+1=2j+1>t$, and $D(x_1)=-t|\theta|=t\theta$.
\item $x_4$: $\fr{-qx_4}=1-2|\theta|$ and $m(G)\equiv p^{-1}-2\equiv q'-2$.
So $\mu=\alpha'$ and $D(x_4)=-2t|\theta|=2t\theta$.
\item $x_2$: $qx_2=\tfrac{3q-p}2+\tfrac{|\theta|}2$, so $\fr{-qx_2}=1-\tfrac{|\theta|}2$ and $2m(G)\equiv(2-p)p^{-1}\equiv2q'-1$.
As before $\mu\in\{\alpha'+\tfrac12,\alpha'-\tfrac12\}$, both larger than $2j$, so $D(x_2)=-\tfrac{t|\theta|}2=\tfrac{t\theta}2$.
\item $x_3$: $qx_3=\bigl(2q-\tfrac{p+1}2\bigr)+\tfrac{1+|\theta|}2$, so $G=2q-\tfrac{p+1}2+1$ and $\fr{-qx_3}=\tfrac{1-|\theta|}2$.
Then $2m(G)\equiv(1-p)p^{-1}\equiv q'-1$, which is even, so $m(G)=\tfrac{q'-1}2\ge1$.
Hence $\mu=\tfrac{\alpha'}2\in(j+\tfrac14,j+\tfrac13)$, $\fl\mu+1=j+1$, and $D(x_3)=t\tfrac{1-|\theta|}2-\min(t,j+1)$.
\item $x_5$: $qx_5=\bigl(\tfrac{7q-1}2-p\bigr)+\tfrac12+|\theta|$, so $\fr{-qx_5}=\tfrac12-|\theta|$ and $2m(G)\equiv(1-2p)p^{-1}\equiv q'-2$.
Since $q'-2$ is odd, $m(G)=\tfrac{q+q'}2-1$.
Hence $\mu=\tfrac{\alpha'-1}2\in(j-\tfrac14,j-\tfrac16)$, $\fl\mu+1=j$, and $D(x_5)=t(\tfrac12-|\theta|)-\min(t,j)$.
\end{itemize}
In every case the stated value of $m(G)$ lies in $\{1,\dots,q\}$, by (E5) and $4\le q'<q$.
Substituting into~\eqref{eq:e-decomp}, all terms containing $\theta$ cancel, and we obtain the stated formulas.
For instance, for $j$ even,
$D_1=2\bigl(t\theta-\tfrac{t\theta}2-\min(t,j)+\tfrac t2-\tfrac{t\theta}2\bigr)=t-2\min(t,j)$.
\end{proof}

\begin{remark}
For $j=1$ the formulas fail (for example $D_1(2q_1)=D_1(2)=0$), because some residues $m(G)$ above degenerate.
Theorem~\ref{thm:e} and Lemma~\ref{lem:block} were checked numerically: the formulas for $2\le j\le16$ and all $t$ with $400$-digit arithmetic, and by brute force with exact floors at all $N=tq_{3j-2}$ with $2\le j\le6$ (up to $N\approx4.8\cdot10^6$).
The values $D_k(tq_i)$ are not the global extremes, and for $k=2$ they are of smaller order.
Since $\log q_{3j-2}\sim j\log j$, the values in Theorem~\ref{thm:e} are $O(\log N/\log\log N)$, whereas by~\cite{rs,ao} and Proposition~\ref{prop:even} the maximum and the minimum of $D_2(M)$ over $M\le N$ are of order $(\log N/\log\log N)^2$.
For example, $\max_{N\le10^9}D_1(N)=23$ and $\min_{N\le10^9}D_1(N)=-18$, whereas Theorem~\ref{thm:e} gives $D_1(83\,128\,184)=8$.
\end{remark}

\section{Quadratic \texorpdfstring{$r$}{r}: automata and block sums}\label{sec:quad}

\subsection{Ostrowski numeration}
Let $\gamma=\fr r=[0;a_1,a_2,\dots]$.
Every $n\ge0$ has a unique Ostrowski representation~\cite{ostrowski,berthe,rockett-szusz}
\[
  n=\sum_{i\ge0}b_iq_i,\qquad 0\le b_0\le a_1-1,\quad 0\le b_i\le a_{i+1},\quad b_i=a_{i+1}\Rightarrow b_{i-1}=0\ (i\ge1),
\]
and for words of a fixed length the lexicographic order of representations is the numerical order.

\begin{lemma}\label{lem:floor}
Let $n\ge1$ have Ostrowski representation $(b_i)$ and let $j$ be the least index with $b_j\ne0$.
Then $\fl{n\gamma}=\sum_ib_ip_i-\one[j\text{ odd}]$.
\end{lemma}

This is classical: it combines $\fl{(n+1)\gamma}=\sum_ib_ip_i$, which is~\cite[Theorem~3]{brown} and goes back to~\cite{fls}, with the description $\fl{(n+1)\gamma}-\fl{n\gamma}=\one[j\text{ odd}]$ of the characteristic word of $\gamma$, which follows from~\cite[Theorem~2]{brown}.
We include the short proof for completeness.

\begin{proof}
We have $n\gamma-\sum b_ip_i=\sum_{i\ge j}b_i\theta_i$.
Let $T_m=(-1)^m\sum_{i\ge m}b_i\theta_i$, so that $T_m=b_m|\theta_m|-T_{m+1}$.
By downward induction, $-|\theta_m|<T_m<|\theta_{m-1}|$ for all $m\ge0$.
The case of large $m$ is trivial.
For the induction step, $T_m>(b_m-1)|\theta_m|\ge-|\theta_m|$, and by~\eqref{eq:rec} $T_m<b_m|\theta_m|+|\theta_{m+1}|\le|\theta_{m-1}|$ for $m\ge1$.
For $m=0$ we get $T_0<(a_1-1)\gamma+|\theta_1|=1-\gamma$.
Now $T_j>(b_j-1)|\theta_j|\ge0$ and $T_j<1$.
So $\sum_{i\ge j}b_i\theta_i=(-1)^jT_j$ lies in $(0,1)$ if $j$ is even and in $(-1,0)$ if $j$ is odd.
\end{proof}

We use Walnut's version of this numeration~\cite{bss,shallit-book}, which omits the redundant place $q_0=q_1=1$ when $a_1=1$.
The places $P_0<P_1<\cdots$, the digit bounds $\beta(j)$, and the rule ``a digit equal to $\beta(j)$ at position $j\ge1$ is followed by $0$'' are as follows:
\begin{itemize}
\item $r=\varphi$, $\gamma=[0;\overline1]$: $P_j=F_{j+2}$ ($1,2,3,5,\dots$) and $\beta(j)=1$ (Zeckendorf).
\item $r=\sqrt2$, $\gamma=[0;\overline2]$: $P_j=q_j$ ($1,2,5,12,29,\dots$), $\beta(0)=1$ and $\beta(j)=2$ for $j\ge1$.
\item $r=\sqrt3$, $\gamma=[0;\overline{1,2}]$: $P_j=q_{j+1}$ ($1,3,4,11,15,41,\dots$), $\beta(j)=2$ for even $j$ and $\beta(j)=1$ for odd $j$.
\end{itemize}
Let $A_K(n)$ be the number whose representation is that of $n$ followed by $K$ zeros.
From $p_i=q_{i-1}$ ($\varphi$, $\sqrt2$), $p_i=q_i-q_{i-2}$ and $q_{i+2}=4q_i-q_{i-2}$ ($\sqrt3$), Lemma~\ref{lem:floor} becomes
\begin{equation}\label{eq:walnut-floor}
  \fl{n\gamma}=A_K(n)-c\,n-\one[\mathrm{odd}(n)],\qquad (K,c)=(1,1),\ (1,2),\ (2,3)
\end{equation}
for $r=\varphi,\sqrt2,\sqrt3$.
Here $\mathrm{odd}(n)$ means that the lowest nonzero digit of $n$ has odd index $j$ in the notation of Lemma~\ref{lem:floor}.
In Walnut's positions this is an even position for $\varphi$ and $\sqrt3$, and an odd position for $\sqrt2$.

Expressing $\fl{n\gamma}$ in Ostrowski numeration by appending zeros, and using this in Walnut, is due to Schaeffer, Shallit and Zorcic~\cite{ssz}.
If $\gamma=[0;\overline{a_1,\dots,a_K}]$ is purely periodic, then $A_K(n-1)=q_K(n-1)+q_{K-1}\fl{n\gamma}$ for $n\ge1$~\cite[Theorem~5]{ssz}.
In our three cases $q_{K-1}=1$ and $q_K=c$, so this reads $\fl{n\gamma}=A_K(n-1)-c\,(n-1)$.
Applied to $n+1$ it gives $\fl{(n+1)\gamma}=A_K(n)-c\,n$, and~\eqref{eq:walnut-floor} follows from this and $\fl{(n+1)\gamma}-\fl{n\gamma}=\one[\mathrm{odd}(n)]$.

\subsection{Automata}
For each of the nine pairs $(r,k)$ with $r$ quadratic in Table~\ref{tab:words} we constructed an msd-first DFAO $M_{r,k}$ over Walnut's Ostrowski numeration.
The construction comes from a comparison automaton for the three cut points of Lemma~\ref{lem:rot}, followed by Moore minimization and, for $r=\varphi,\sqrt3$, a shift of the outputs that accounts for the place $q_0$ omitted by Walnut; the script that builds the automata is in the code archive.
Table~\ref{tab:quad} gives the number of states of the automata $M_{r,k}$ as used.
Four of them are not minimal as DFAOs over the full digit alphabet: the minimal DFAOs for $(\varphi,2)$, $(\varphi,3)$, $(\varphi,4)$ and $(\sqrt3,2)$ have $8$, $11$, $13$ and $18$ states, and the other five are minimal.
Nothing below depends on minimality.

\begin{proposition}\label{prop:walnut}
For each of the nine pairs, $M_{r,k}(n)=s_k(n)$ for all $n\ge0$.
\end{proposition}

\begin{proof}
We work in Walnut~\cite{walnut}.
Let \texttt{app} be the two-track automaton that accepts exactly the empty pair and the pairs of words $(0^Kw,\,w0^K)$, where $w$ is any word over the digit alphabet.
As a relation between numbers, for $m\ge1$ (the only case in which it is used below), $\mathtt{app}(m,l)$ holds if and only if $l=A_K(m)$.
Let \texttt{odd} be the regular expression for $\mathrm{odd}(n)$.
Then~\eqref{eq:walnut-floor} defines a first-order formula for $y=\fl{m\gamma}$, as in~\cite{ssz}.
For $r=\sqrt2$ it reads
\begin{verbatim}
def fl_rt2 "?msd_rt2 (m=0 & y=0) | (m>=1 & El $app_rt2(m,l) &
   (($odd_rt2(m) & y+1+2*m=l) | (~$odd_rt2(m) & y+2*m=l)))";
\end{verbatim}
and then, for example for $k=2$,
\begin{verbatim}
def s_rt2_2 "?msd_rt2 Ea,b,c,x,y x=2*n+2 & y=n+2 & $fl_rt2(x,a)
   & $fl_rt2(y,b) & $fl_rt2(n,c) & a=b+c+1";
eval s_rt2_2_dfao "?msd_rt2 An (Srt22[n]=@1 <=> $s_rt2_2(n))";
\end{verbatim}
Walnut returns \texttt{TRUE} for all nine statements of the last form, and for the nine statements asserting that, for every $n$, the same formula gives $a\in\{b+c,\,b+c+1\}$, that is, $s_k(n)\in\{0,1\}$.
It also returns \texttt{TRUE} for the checks that each \texttt{fl} is a function, is nondecreasing with steps $0$ or $1$, and has the right values at ten sample points; this makes $27$ statements in all.
The nine negated control statements return \texttt{FALSE}.
The automata \texttt{app} were checked independently against their defining relation by an exhaustive product search.
\end{proof}

\begin{table}[ht]
\centering
\begin{tabular}{c|c|ccc}
state & output & on $0$ & on $1$ & on $2$\\ \hline
0 & 0 & 0 & 1 & 2\\
1 & 0 & 3 & 1 & 2\\
2 & 0 & 4 & 5 & 5\\
3 & 1 & 1 & 1 & 2\\
4 & 1 & 6 & 7 & 2\\
5 & 0 & 5 & 5 & 5\\
6 & 1 & 3 & 1 & 2\\
7 & 0 & 4 & 1 & 2
\end{tabular}
\caption{The DFAO $M_{\sqrt2,2}$ (msd first, initial state $0$).
For example $4=[20]$ gives the path $0\to2\to4$ and output $s_2(4)=1$.
State $5$ is reached only by invalid words.}
\label{tab:dfao}
\end{table}

\subsection{Block sums}
Fix one of the nine cases, with DFAO $M=(Q,\delta,q_0,\tau)$, and put $g(q)=2\tau(q)-1$.
Let $\mathcal L_i$ be the set of valid words of length $i$ occupying positions $i-1,\dots,0$.
For $q\in Q$ and $c\in\{0,1\}$ let
\[
  W_i(q,c)=\sum_x g\bigl(\delta^*(q,x)\bigr),
\]
where $x$ runs over $\mathcal L_i$, and additionally $x$ must start with $0$ if $c=1$.
Here $W_0(q,c)=g(q)$, and
\begin{equation}\label{eq:Wrec}
  W_i(q,c)=\sum_{e=0}^{c\,?\,0\,:\,\beta(i-1)}W_{i-1}\bigl(\delta(q,e),\one[e=\beta(i-1)]\bigr).
\end{equation}

\begin{lemma}\label{lem:dp}
Let $N\ge1$ have a valid word $N_{L-1}\cdots N_0$ (leading zeros allowed), and let $\sigma_i=\delta^*(q_0,N_{L-1}\cdots N_{i+1})$.
Then
\[
  D_k(N-1)=1+\sum_{i=0}^{L-1}\ \sum_{d=0}^{N_i-1}W_i\bigl(\delta(\sigma_i,d),0\bigr).
\]
\end{lemma}

\begin{proof}
The integers $0\le n<N$ correspond to the valid words of length $L$ that are lexicographically smaller than $N_{L-1}\cdots N_0$.
Such a word agrees with $N$ above some position $i$ and has a digit $d<N_i$ at position $i$.
Since $d<N_i\le\beta(i)$, the digit $d$ imposes no restriction on the lower digits, which form an arbitrary word of $\mathcal L_i$.
(If $N_{i+1}=\beta(i+1)$, then $N_i=0$ and position $i$ contributes nothing.)
This gives $\sum_{n=0}^{N-1}(2s_k(n)-1)$, and the term $n=0$ equals $-1$.
\end{proof}

For $i\ge0$ let $S_i$ be the set of pairs $\bigl(\delta^*(q_0,y),\one[y_i=\beta(i)]\bigr)$, where $y$ runs over the valid words occupying positions $L-1,\dots,i$ for some $L>i$, and $y_i$ is the digit of $y$ at position $i$.
Lemma~\ref{lem:dp} evaluates $W_i$ only on $S_i$, and~\eqref{eq:Wrec} maps $S_i$ into $S_{i-1}$, so only the restrictions $W_i|_{S_i}$ matter.
Since the bounds $\beta(j)$, $j\ge1$, have period $p$, we have $S_{i+p}=S_i$ for $i\ge1$.
(For $r=\sqrt3$, where $p=2$, the sets $S_i$ for even and odd $i$ differ, and on their union the values $W_i$ are unbounded; the restriction to $S_i$ is essential.)
On these sets the vectors $W_i$ are eventually periodic.
This is a finite check: for $i\ge2$ the coefficients in~\eqref{eq:Wrec} depend only on $i$ modulo $p$.
Hence, if $W_{i_0+T}=W_{i_0}$ on $S_{i_0}$ for some $i_0\ge2$ and some multiple $T$ of $p$, then $W_{i+T}=W_i$ on $S_i$ for all $i\ge i_0$.
Table~\ref{tab:quad} lists $i_0$ and $T$.
In all nine cases $|W_i(q,c)|\le3$ for $(q,c)\in S_i$.

\begin{lemma}[insertion]\label{lem:insert}
Let $u,w,v$ be words such that $p$ and $T$ divide $|w|$, $|w|+|v|\ge i_0$, $\delta^*(q_0,uw)=\delta^*(q_0,u)$, and $uwv$ and $uwwv$ are valid.
Let $N_m$ be the integer with representation $uw^mv$, and assume that $N_1\ge1$, that is, that $uwv$ has a nonzero digit.
Then $D_k(N_{m}-1)=D_k(N_1-1)+(m-1)C$ for all $m\ge1$, where $C=D_k(N_2-1)-D_k(N_1-1)$.
\end{lemma}

\begin{proof}
Validity is a local condition on adjacent digits with $p$-periodic bounds at positions $j\ge1$.
So $uw^mv$ is valid for all $m\ge1$.
Like $uwv$, it has a nonzero digit, so $N_m\ge1$, and Lemma~\ref{lem:dp} applies to $N_m$.
Compare the sums of Lemma~\ref{lem:dp} for $N_{m+1}$ and $N_m$.
The part coming from the final $w^mv$ is identical, because it is read from the same state $\delta^*(q_0,uw)=\delta^*(q_0,u)$ at the same positions.
The part coming from $u$ is read from the same states, at positions shifted by $|w|$, all of them at least $|wv|\ge i_0$.
By periodicity it is also identical.
What remains is the contribution of the inserted copy of $w$.
It is read from the state $\delta^*(q_0,u)$ at positions whose residues modulo $T$ do not depend on $m$, all at least $i_0$.
So it is a constant $C$.
\end{proof}

\begin{table}[ht]
\centering
\small
\begin{tabular}{cc|cc|ll|ll}
\toprule
$r$ & $k$ & states & $(i_0,T)$ & up: $u,w,v$ & $D_k(N_m-1)$ & down: $u,w,v$ & $D_k(N_m-1)$\\
\midrule
$\varphi$ & 1 & 8 & (2,6) & 1, 000101, --- & $2m$ & 1, 000001, 0 & $-1-2m$\\
$\varphi$ & 2 & 9 & (2,6) & 1, 000001, --- & $m$ & 1, 000101, 0 & $1-m$\\
$\varphi$ & 3 & 13 & (2,6) & 1, 000001, --- & $m$ & 1, 000001, 0 & $1-2m$\\
$\varphi$ & 4 & 15 & (3,6) & 1, 000001, --- & $m$ & 1, 000001, 000 & $-m$\\
$\sqrt2$ & 1 & 8 & (2,2) & 1, 11, --- & $m$ & 1, 01, --- & $-m$\\
$\sqrt2$ & 2 & 8 & (2,2) & \multicolumn{2}{l|}{$\sup D_2=0$} & 1, 01, --- & $-m$\\
$\sqrt2$ & 4 & 12 & (3,2) & \multicolumn{2}{l|}{$\sup D_4=1$} & 1, 01, 00 & $-m$\\
$\sqrt3$ & 1 & 20 & (2,4) & 1, 0101, --- & $3m$ & 1, 0001, 0 & $-m$\\
$\sqrt3$ & 2 & 20 & (2,4) & 1, 0001, --- & $m$ & 1, 0001, 0 & $-m$\\
\bottomrule
\end{tabular}
\caption{DFAO sizes, periods of the block sums, and digit families $N_m=[uw^mv]$ (Walnut's Ostrowski digits, most significant first; --- is the empty word) with $D_k(N_m-1)$ for all $m\ge1$.}
\label{tab:quad}
\end{table}

\begin{proposition}\label{prop:families}
The families in Table~\ref{tab:quad} satisfy the hypotheses of Lemma~\ref{lem:insert}, and $D_k(N_m-1)$ is as stated for all $m\ge1$.
In particular $D_k$ is unbounded above in the seven cases with an ``up'' family, and unbounded below in all nine cases.
\end{proposition}

\begin{proof}
The hypotheses are finite checks (in particular $N_1\ge1$, since $u=1$), and $D_k(N_1-1)$ and $D_k(N_2-1)$ are computed with Lemma~\ref{lem:dp}; then Lemma~\ref{lem:insert} applies.
We also confirmed the values by brute force with exact integer floors for every $N_m\le2\cdot10^6$.
\end{proof}

For $r=\sqrt2$ the families have closed forms in terms of the Pell numbers $\mathcal P_0=0$, $\mathcal P_1=1$, $\mathcal P_{n+1}=2\mathcal P_n+\mathcal P_{n-1}$ (so that Walnut's places are $P_j=\mathcal P_{j+1}$) and $H_n=\mathcal P_n+\mathcal P_{n-1}$.
Namely $[1(01)^m]=\mathcal P_{2m+2}/2$, $[1(01)^m00]=\mathcal P_{2m+4}/2-1$ and $[1(11)^m]=(H_{2m+2}-1)/2$.
Hence, for all $m\ge1$,
\[
  D_1\Bigl(\tfrac{\mathcal P_{2m+2}}2-1\Bigr)=D_2\Bigl(\tfrac{\mathcal P_{2m+2}}2-1\Bigr)=-m,\qquad D_4\Bigl(\tfrac{\mathcal P_{2m+4}}2-2\Bigr)=-m,\qquad D_1\Bigl(\tfrac{H_{2m+2}-3}{2}\Bigr)=m .
\]
The numbers $\mathcal P_{2m}/2$ are the half-Pell numbers at which, by a theorem of Bruin and Fokkink~\cite[Theorem~5]{bf}, the records of $2D_N(\sqrt2,\tfrac12)$ occur.

\subsection{Exact suprema and growth rates}
For $i\ge0$, $q\in Q$ and $c\in\{0,1\}$, let $V_i(q,c)$ be the largest value, over the words $x\in\mathcal L_i$ allowed by $c$, of the contribution of positions $i-1,\dots,0$ to the sum in Lemma~\ref{lem:dp} when these digits are read from state $q$.
Then $V_0=0$, and
\[
  V_i(q,c)=\max_{e}\Bigl(\sum_{d<e}W_{i-1}\bigl(\delta(q,d),0\bigr)+V_{i-1}\bigl(\delta(q,e),\one[e=\beta(i-1)]\bigr)\Bigr),
\]
with $e$ as in~\eqref{eq:Wrec}.
If the leading digit of $N$ is $e\ge1$ at position $i$, then
\[
  D_k(N-1)\le1+\sum_{d<e}W_i\bigl(\delta(q_0,d),0\bigr)+V_i\bigl(\delta(q_0,e),\one[e=\beta(i)]\bigr),
\]
with equality for a suitable $N$.
This recursion is max-plus linear; at $(q,c)\in S_i$ it involves only $V_{i-1}$ and $W_{i-1}$ on $S_{i-1}$, so its coefficients are periodic in $i$ for $i\ge i_0$.
Therefore, if $V_{i_1+T'}=V_{i_1}+hT'$ on $S_{i_1}$ (the same shift in every coordinate) for some $i_1\ge i_0$ and some multiple $T'$ of $T$, then $V_{i+T'}=V_i+hT'$ on $S_i$ for all $i\ge i_1$.
The same applies to minima.

\begin{proposition}\label{prop:sup}
For $r=\sqrt2$ we have $\sup_{N\ge0}D_2(N)=0$ and $\sup_{N\ge0}D_4(N)=1$.
\end{proposition}

\begin{proof}
For $k=2$ one finds $V_{i+2}=V_i$ for $i\ge2$; for $k=4$, $V_{i+2}=V_i$ for $i\ge3$.
So the supremum is a maximum over the first few positions.
It equals $0$ (at $N=1$, that is, $D_2(0)=0$; also, for example, $D_2(2)=0$) and $1$ (at $D_4(3)=1$), respectively.
\end{proof}

\begin{proposition}\label{prop:growth}
In the notation of Table~\ref{tab:quad}, let $\lambda_\pm$ be the rates below and let $\eta$ be $\log\varphi$, $\log(1+\sqrt2)$ or $\tfrac12\log(2+\sqrt3)$.
Then
\[
  \max_{0\le M\le X}D_k(M)=\frac{\lambda_+}{\eta}\log X+O(1),\qquad \min_{0\le M\le X}D_k(M)=-\frac{\lambda_-}{\eta}\log X+O(1),
\]
whenever the corresponding side is unbounded.
The rates per digit are:
\begin{center}
\begin{tabular}{c|cccc|ccc|cc}
 & $\varphi,1$ & $\varphi,2$ & $\varphi,3$ & $\varphi,4$ & $\sqrt2,1$ & $\sqrt2,2$ & $\sqrt2,4$ & $\sqrt3,1$ & $\sqrt3,2$\\ \hline
$\lambda_+$ & $1/3$ & $1/6$ & $1/3$ & $1/6$ & $1$ & --- & --- & $3/4$ & $1/4$\\
$\lambda_-$ & $1/3$ & $1/6$ & $1/3$ & $1/6$ & $1/2$ & $1/2$ & $1/2$ & $1/4$ & $1/4$
\end{tabular}
\end{center}
\end{proposition}

\begin{proof}
In each case the value vectors reach the linear regime $V_{i+T'}=V_i+hT'$ (and the analogous one for minima), with $h=\pm\lambda_\pm$.
Hence $\max_{1\le N<P_L}D_k(N-1)=\lambda_+L+O(1)$.
Since $\log P_L=\eta L+O(1)$, the claim follows.
\end{proof}

\begin{remark}\label{rem:rates}
For even $k$ the rates in Proposition~\ref{prop:growth} also follow from known results on $D_N(r,\tfrac12)$, since $D_k(N)=-2D_N(r,\tfrac12)+O(1)$ by Proposition~\ref{prop:even}.
The formulas of Ro\c{c}adas and Schoi\ss{}engeier~\cite{rs} give all five even-$k$ columns.
For $r=\varphi$ and $r=\sqrt3$ the least $N$ at which $2D_N(r,\tfrac12)=S_N(2r)$ takes a given value is given exactly by~\cite[Corollary~4.5]{ors}.
Moreover,~\cite[Theorem~3.1]{ors} with $m=2$ (which applies since $2\varphi=(\sqrt5-1)+2$ and $S_N(\beta+2)=S_N(\beta)$) gives $|S_N(2\varphi)|\le\log N/(2\log(2+\sqrt5))+1$, and $1/(2\log(2+\sqrt5))=1/(6\log\varphi)$ is the constant $\lambda_\pm/\eta$ for $(\varphi,2)$ and $(\varphi,4)$.
For $r=\sqrt2$ the records of $S_N(2\sqrt2)$ at the half-Pell numbers~\cite[Theorem~5]{bf} give $\lambda_-=1/2$ for $(\sqrt2,2)$ and $(\sqrt2,4)$.
The rates for odd $k$ and the exact suprema in Proposition~\ref{prop:sup} do not follow from these results.
\end{remark}

For example, $\lambda_+/\eta=1/\log(1+\sqrt2)\approx1.135$ for $(\sqrt2,1)$, and indeed $\max_{N\le10^9}D_1(N)=24$ while $1.135\log10^9\approx23.5$.
The families of Table~\ref{tab:quad} were chosen to be short, not optimal: for $(\varphi,3)$ and $(\sqrt2,1)$ the ``up'' family gains $1/6$ and $1/2$ per digit, less than $\lambda_+$.

\section{Results for the 22 sequences}\label{sec:results}

Theorem~\ref{thm:main} follows from Lemma~\ref{lem:equiv}, Theorem~\ref{thm:e} and Propositions~\ref{prop:families} and~\ref{prop:sup}.
Table~\ref{tab:results} lists the verdict for each sequence.
For the three false conjectures and for A327224, the bounds follow from Lemma~\ref{lem:equiv} and Proposition~\ref{prop:sup}.
For $k=2$: $2m-u(m)=D_2(u(m))\le0$, and $2m-v(m)=-D_2(v(m))=1-D_2(v(m)-1)\ge1$.
For $k=4$: $2m-u(m)\le1$, and $2m-v(m)=1-D_4(v(m)-1)\ge0$.
The other side is unbounded by Proposition~\ref{prop:families}.

\begin{table}[ht]
\centering
\begin{tabular}{llccl}
\toprule
sequence & $r$ & $k$ & letter & verdict\\
\midrule
A327175 & $\varphi$ & 1 & 0 & true (Prop.~\ref{prop:families})\\
A327176 & $\varphi$ & 1 & 1 & true (Prop.~\ref{prop:families})\\
A327178 & $\sqrt2$ & 1 & 0 & true (Prop.~\ref{prop:families})\\
A327179 & $\sqrt2$ & 1 & 1 & true (Prop.~\ref{prop:families})\\
A327181 & $\sqrt3$ & 1 & 0 & true (Prop.~\ref{prop:families})\\
A327182 & $\sqrt3$ & 1 & 1 & true (Prop.~\ref{prop:families})\\
A327184 & $e$ & 1 & 0 & true (Thm.~\ref{thm:e})\\
A327185 & $e$ & 1 & 1 & true (Thm.~\ref{thm:e})\\
A327203 & $\varphi$ & 2 & 0 & true (Cor.~\ref{cor:yz}, Prop.~\ref{prop:families})\\
A327204 & $\varphi$ & 2 & 1 & true (Cor.~\ref{cor:yz}, Prop.~\ref{prop:families})\\
A327206 & $\sqrt2$ & 2 & 0 & \textbf{false}: $2n-a(n)\ge1$, equality at $n=1$\\
A327207 & $\sqrt2$ & 2 & 1 & \textbf{false}: $2n-a(n)\le0$, i.e.\ $a(n)\ge2n$, equality at $n=1$\\
A327209 & $\sqrt3$ & 2 & 0 & true (Cor.~\ref{cor:yz}, Prop.~\ref{prop:families})\\
A327210 & $\sqrt3$ & 2 & 1 & true (Cor.~\ref{cor:yz}, Prop.~\ref{prop:families})\\
A327212 & $e$ & 2 & 0 & true (Cor.~\ref{cor:yz}, Thm.~\ref{thm:e})\\
A327213 & $e$ & 2 & 1 & true (Cor.~\ref{cor:yz}, Thm.~\ref{thm:e})\\
A327217 & $\varphi$ & 3 & 0 & true (Prop.~\ref{prop:families})\\
A327218 & $\varphi$ & 3 & 1 & true (Prop.~\ref{prop:families})\\
A327220 & $\varphi$ & 4 & 0 & true (Cor.~\ref{cor:yz}, Prop.~\ref{prop:families})\\
A327221 & $\varphi$ & 4 & 1 & true (Cor.~\ref{cor:yz}, Prop.~\ref{prop:families})\\
A327223 & $\sqrt2$ & 4 & 0 & \textbf{false}: $2n-a(n)\ge0$, equality at $n=2$\\
A327224 & $\sqrt2$ & 4 & 1 & true: $a(n)-2n$ is unbounded above, and $a(n)-2n\ge-1$\\
\bottomrule
\end{tabular}
\caption{Verdicts. ``Letter'' is the letter whose positions the sequence lists.}
\label{tab:results}
\end{table}

The bounds are attained many times.
Among the terms with $a(n)\le2\cdot10^6$, equality holds for $8747$ indices $n$ in A327206, for $8746$ in A327207, and for $4373$ each in A327223 and A327224.
In the data listed in the entries ($63$, $62$, $63$ and $62$ terms), the observed ranges of $2n-a(n)$ for A327206, A327207, A327223 and A327224 are $[1,3]$, $[-2,0]$, $[0,3]$ and $[-2,1]$, so the one-sidedness is already visible there.

\begin{remark}[OEIS comments]\label{rem:oeis}
Each position sequence carries a second comment, ``Positive integers $k$ such that A\dots$(k)=0$'' or ``$=1$''.
These comments were added by Jianing Song.
For A327182, A327204 and A327207, which list positions of $1$s, the comment says ``$=0$'' instead of ``$=1$''.
\end{remark}

\FloatBarrier
\section{Verification}\label{sec:verification}

All computations are in the accompanying code archive.
\begin{itemize}
\item \texttt{exact.c}: brute-force $D_k(N)$ for $N\le10^9$ with exact floors. For $\varphi$, $\sqrt2$ and $\sqrt3$ the floors are computed with integer square roots. For $e$, $me$ is computed in quadruple precision with error below $10^{-23}$. Since $q_{22}=312\,129\,649\le2\cdot10^9+4<q_{23}=5\,155\,334\,720$ (denominators of the convergents of $e$), every argument $1\le m\le2\cdot10^9+4$ satisfies $\|me\|\ge|q_{22}e-p_{22}|>10^{-10}$, where $\|x\|$ is the distance from $x$ to the nearest integer, so these floors are exact as well.
It covers all four $r$ and $k\le4$.
The extremes agree with Theorem~\ref{thm:main}; for instance $\max_{N\le10^9}D_2(N)=0$ and $\max_{N\le10^9}D_4(N)=1$ for $r=\sqrt2$.
\item Walnut command files \texttt{fk\_gold.txt}, \texttt{fk\_rt2.txt}, \texttt{fk\_rt3.txt}, the automata \texttt{app\_*.txt}, and the DFAOs \texttt{S*.txt}.
All $27$ statements return \texttt{TRUE} and the $9$ controls return \texttt{FALSE}, re-run from a clean installation of the development version 8.0-alpha of Walnut (commit 49a16a3, September 2026).
The script \texttt{construct/walnut\_gen.py} rebuilds these files and the DFAOs.
\item \texttt{check\_quadratic.py}: an independent checker for Section~\ref{sec:quad}.
It uses only the DFAO files, exact integer floors and its own implementation of the numeration.
It checks that the valid words of length $L\le12$ are exactly the greedy representations of $0,\dots,P_L-1$, in order.
It checks the DFAOs against $s_k(n)$ for $n<10^6$, and the formula of Lemma~\ref{lem:dp} against brute force for $N<6\cdot10^4$.
It then verifies the periodicity of $W$, searches for the families, verifies the hypotheses of Lemma~\ref{lem:insert}, and computes the value recursions (Propositions~\ref{prop:sup} and~\ref{prop:growth}).
\item \texttt{check\_append.py}: product check of the \texttt{app} automata.
\item \texttt{e\_verify.py} with \texttt{blockformula.py}: Lemma~\ref{lem:block} on $738$ random instances (six values of $\alpha$), and Theorem~\ref{thm:e} as described above.
\item \texttt{check\_reduction.py}: Proposition~\ref{prop:even} numerically for $N\le10^6$.
\item \texttt{check\_oeis.py}: the OEIS data of all $33$ entries agree with~\eqref{eq:s}.
\item \texttt{attain.py}: the attainment counts in Section~\ref{sec:results}.
\end{itemize}

\section*{Acknowledgements}
The conjectures treated here were posed by Clark Kimberling.
The author used an AI assistant (Claude, Anthropic) for literature searches, computations, and drafting, and checked the proofs and computations.

\end{document}